\documentclass[journal]{IEEEtran}
\ifCLASSINFOpdf
\else
\fi

\usepackage{amsmath}
\usepackage{subfig}
\usepackage{graphicx}
\usepackage{array}
\newtheorem{theorem}{Theorem}
\begin{document}
%
\title{A fast and effective method for a Poisson denoising model with total variation}
%
%
%

\author{Wei~Wang
	and~Chuanjiang~He
\thanks{The authors are with the College of Mathematics and Statistics, Chongqing University, Chongqing 401331, China.}
\thanks{E-mail: 419714480@qq.com; cjhe@cqu.edu.cn}
\thanks{Manuscript received ; revised }}
%
%

\markboth{IEEE SIGNAL PROCESSING LETTERS, VOL. N, NO. N, June 2016}{Shell \MakeLowercase{\textit{et al.}}: Bare Demo of IEEEtran.cls for IEEE Journals}
%



\maketitle

\begin{abstract}
In this paper, we present a fast and effective method for solving the Poisson-modified total variation model proposed in [9]. The existence and uniqueness of the model are again proved using different method. A semi-implicit difference scheme is designed to discretize the derived gradient descent flow with a large time step and can guarantee the restored image to be strictly positive in the image domain. Experimental results show the efficiency and effectiveness of our method.
\end{abstract}


%
\IEEEpeerreviewmaketitle

\section{Introduction}
%
%
%
%
\IEEEPARstart{P}{oisson}  noise, also known as photon noise, is a basic form of uncertainty associated with the measurement of light. An image sensor measures scene irradiance by counting the number of photons incident on the sensor over a given time interval. The photon counting is a classic Poisson process that follows Poisson distribution [1]. 
Poisson noise removal is a fundamental task for many imaging applications where images are generated by photon-counting devices such as computed tomography (CT), magnetic resonance imaging (MRI) and astronomical imaging. Many methods and algorithms have been proposed for Poisson denoising [2-21]. In this paper, we focus on solving the variational Poisson denoising model proposed in [9].

In [9], along the lines of the famous ROF model [22], Le, Chartrand and Asaki proposed the following Poisson denoising model with total variation regularization (called LCA model in this paper). In detail, if $f = f(x,y)$  ( $(x,y) \in \Omega $, a bounded, open subset of  ${R^2}$) is an original image with Poisson noise, then the reconstructed image ${u^ * }$  is obtained by

\begin{equation}\label{key}
{u^ * } = \arg \mathop {\inf }\limits_u E(u) = \int_\Omega  {\left| {\nabla u} \right|}  + \beta \int_\Omega  {(u - f\log u)},
\end{equation}
where the functional $E$  is defined on the set of $u \in BV(\Omega )$  such that  $\log u \in {L^1}(\Omega )$; in particular, $u$  must be positive almost everywhere (a.e.) in  $\Omega $.

The authors used gradient descent with the forward-time central-space finite difference scheme to solve problem (1). They implemented a straightforward, discretized version of the following PDE:
\begin{equation}\label{key}
\frac{{\partial u}}{{\partial t}} = div\left( {\frac{{\nabla u}}{{\left| {\nabla u} \right|}}} \right) + \beta \left( {\frac{f}{u} - 1} \right)\; \text{with} \;\frac{{\partial u}}{{\partial \vec n}} = 0\; \text{on}\; \partial \Omega. 
\end{equation}            
Spatial derivatives are computed with standard centered difference approximations. The quantity $\left| {\nabla u} \right|$  is replaced with  $\sqrt {{{\left| {\nabla u} \right|}^2} + \varepsilon } $ for a small, positive  $\varepsilon $. The time evolution is done with fixed time step  $\tau $, until the change in $u$  is sufficiently small.
However, this numerical scheme has two main drawbacks: 1) The reconstructed image $\tilde u$  obtained by this scheme cannot be guaranteed to be positive a.e. in $\Omega $  mathematically; in fact, our experiments show that the reconstructed image $\tilde u$  is sign-changing in  $\Omega $. The functional $E$  is singular in the non-positive orthant of the sign-changing solution $\tilde u$  due to the presence of $\log u$  in the functional  $E$. Therefore, $\tilde u$  is not the best approximation to ${u^ * }$  defined by (1). 2) The time step $\tau $  must be chosen small enough to ensure the stability of the used explicit numerical scheme due to the CFL condition.

After the LCA model, many algorithms were proposed to solve the LCA model. In [10], Chan and Chen proposed a multilevel algorithm for efficiently solving the LCA model. However, this method also confronts the problem that the reconstructed image is sign-changing in  $\Omega $ and it spends a little long time. In [13], Figueiredo and Bioucas-Dias used an alternating direction method of multipliers to solve the LCA model. To address the problem that the reconstructed image is sign-changing in  $\Omega $, they replaced $u$  with the projection $\max (u,0)$  during iteration. However, $E$  is still meaningless because of $\max (u,0)=0$ for $u\le0$.

In this letter, we first prove the existence and uniqueness of the solution for the LCA model by a different method. Then  a semi-implicit difference scheme is designed to solve numerically the derived gradient descent flow of the LCA model. The proposed scheme can guarantee the restored image to be strictly positive in the image domain and is stable for a large time step.

\section{The proposed algorithm}
In this paper, we propose a new method to solve the  Poisson denoising model [9]:

\begin{equation}\label{key}
\mathop {\inf }\limits_{u \in G(\Omega )} E(u) = \int_\Omega  {\left| {\nabla u} \right|}  + \beta \int_\Omega  {(u - f\log u)},
\end{equation}                    
where $G(\Omega ) = \left\{ {u \in BV(\Omega ):u > 0,\;a.e.\;{\rm{in}}\;\Omega } \right\}$  is a subset of  $BV(\Omega )$ and $f$ is an original image with Poisson noise.

Inspired by the proof of Theorem 4.1 in [23], in the following we give a new proof of the existence and uniqueness of solution for problem (3), which differs from the proof in [9].
\begin{theorem}
	If  $f \in G(\Omega )$, then problem (1) has exactly one solution.
\end{theorem}
\begin{IEEEproof}
	Let $M = \sup f$ and  $m = \inf f$. Since $H(u) = u - f\log u$  decreases if $u \in (0,f)$  and increases if  $u \in (f, + \infty )$, we have 
	\begin{align*}
	&\int_\Omega  {\min (u,M) - f\log (\min (u,M))} \\
	&= \left( {\int_{\Omega (u \le M)} { + \int_{\Omega (u > M)} {} } } \right)\min (u,M) - f\log (\min (u,M))\\
	&\le \int_\Omega  {u - f\log u} 
	\end{align*}
	and similarly,
	$$\int_\Omega  {\max (u,m) - f\log (\max (u,m))}  \le \int_\Omega  {u - f\log u}.$$
	Since
	\begin{align*}
	\int_\Omega  {\left| {\nabla \min (u,M)} \right|}  &= \left( {\int_{\Omega (u \le M)} { + \int_{\Omega (u > M)} {} } } \right)\left| {\nabla \min (u,M)} \right|\\
	&= \int_{\Omega (u \le M)} {\left| {\nabla u} \right|}  \\
	&\le \int_\Omega  {\left| {\nabla u} \right|} 
	\end{align*}
	and  
	$$\int_\Omega  {\left| {\nabla \max (u,m)} \right|}  \le \int_\Omega  {\left| {\nabla u} \right|}, $$ 
	we have $$E(\min (u,M)) \le E(u)$$
	and  $$E(\max (u,m)) \le E(u).$$ Thus we can assume 
	$0 \le m \le u \le M$ a.e. in $\Omega $.
	
	Since  $f \in G(\Omega )$, we have $\log f \in {L^1}(\Omega )$  and $E(u) \ge \beta \int_\Omega  {(u - f\log u)}  \ge \beta \int_\Omega  {(f - f\log f)} $,
	which implies that $E$ is bounded below in  $G(\Omega )$. Let $\left\{ {{u_n}} \right\}$  be the minimization sequence of problem (3) such that 
	$$\mathop {\lim }\limits_{n \to \infty } E({u_n}) = \mathop {\inf }\limits_{u \in G(\Omega )} E(u): = {E_0}.$$
	Then, there is an $N$ such that, for every  $n > N$,
	$$\int_\Omega  {\left| {\nabla {u_n}} \right|}  + \beta \int_\Omega  {({u_n} - f\log {u_n}} ) \le {E_0} + 1,$$
	which implies 
	\begin{align*}
	\int_\Omega  {\left| {\nabla {u_n}} \right|}  &\le {E_0} + 1 - \beta \int_\Omega  {({u_n} - f\log {u_n}} )\\
	&\le  {E_0} + 1 - \beta \int_\Omega  {f(1 - \log f} )
	\end{align*}
	Recalling that  $M \le {u_n} \le m$, thus $\left\{ {{u_n}} \right\}$  is bounded in $BV(\Omega )$  which implies there exists a ${u^ * } \in BV(\Omega )$  such that, up to a subsequnce, 
	${u_n} \to {u^ * }$  weakly in ${L^2}(\Omega )$  and strongly in  ${L^1}(\Omega )$.
	Since  ${u_n} \ge 0$ a.e. in $\Omega$, we have ${u^ * } \ge 0$ a.e. in $\Omega$. We further have ${u^ * } > 0$ a.e. in $\Omega$; otherwise, up to a sequence,
	$$\mathop {\lim }\limits_{n \to \infty } \int_\Omega  { - f\log ({u_n})}  = \int_\Omega  { - f\log ({u^ * })}  =  + \infty $$
	which contradicts with the fact that  $\mathop {\lim }\limits_{n \to \infty } E({u_n}) = {E_0}.$  Thus  ${u^ * } \in G(\Omega )$.  
	Thanks to the lower semi-continuity of the total variation and Fatou’s lemma, we get that $u^{*}$  is a solution of problem (3). The uniqueness of  the solution is guaranteed by the strict convexity of problem (3).	
\end{IEEEproof}	

In what follows, we design a semi-implicit difference scheme for the discretization of the gradient descent flow for problem (3).

The gradient descend flow of problem (3) is 
\begin{equation}\label{key}
\frac{{\partial u}}{{\partial t}} = div(\frac{{\nabla u}}{{\left| {\nabla u} \right|}}) - \beta (1 - \frac{f}{u})
\end{equation}

To guarantee that the restored image is positive, we design the following semi-implicit difference scheme to discretize equation (4):
\begin{equation}\label{key}
\frac{{{u^{n{\rm{ + 1}}}} - {u^n}}}{\tau } = div\left( {\frac{{\nabla {u^n}}}{{\left| {\nabla {u^n}} \right|}}} \right) - \beta \left( {1 - \frac{f}{{{u^{n + 1}}}}} \right).
\end{equation}
Equation (5) can be rewritten as 
\begin{equation}\label{key}
{({u^{n{\rm{ + 1}}}})^2} - \left( {{u^n} + \tau (div(\frac{{\nabla {u^n}}}{{\left| {\nabla {u^n}} \right|}}) - \beta )} \right){u^{n{\rm{ + 1}}}} - \beta \tau f = 0.
\end{equation}
The restored image is given as the positive solution of equation (6):
\begin{equation}\label{key}
{u^{n + 1}} = \frac{{ - {a_n} + \sqrt {a_n^2 - 4b} }}{2}
\end{equation}
where 
${a_n} =  - {u^n} - \tau (div(\frac{{\nabla {u^n}}}{{\left| {\nabla {u^n}} \right|}}) - \beta )$,  $b =  - \beta \tau f$.


\section{Experiments}
In this section, we present some experimental results to show the performance and the effectiveness of the proposed numerical scheme, in comparison to other relevant numerical schemes in [9], [10] and [13] in terms of quality and time. Ten images are chosen as test images, which are shown in Figure 1. The corresponding Poisson noisy images were generated from the test images by using the Matlab command `imnoise` with noise type parameter `poisson`. Since Poisson noise depends only on the intensity of the image, there is no extra parameter in the Matlab command. The stop criteria for our method is set as  $\left| {(E({u^{n + 1}}) - E({u^n}))/E({u^{n + 1}})} \right|\le tolerance $. 

The parameters for the proposed numerical scheme are set as  $\beta  = 10$,  $\tau  = {\rm{0}}{\rm{.7}}$ and $tolerance = 3.0e-4$ . The parameters for the numerical scheme in [9] are set as  $\beta  = 10$, $\tau  = 0.01$  and  iteration numbers are all set 30 for all experiments. The parameters for the multilevel algorithm [10] are set as  $\alpha = 0.05$, $tol = 1.0e - 3$; by our experiments, $\alpha = 0.05$  is averagely the best value for the ten test images. The parameters for the method [13] are set as the same values as in [13]:  $\tau {\rm{ = 0}}{\rm{.1}}$,  $\mu {\rm{ = 60}}\tau {\rm{/}}M$, where $M$  is the maximal value of the Poisson image. The inner and outer iteration numbers of the method [13] are 10 and 6, respectively, which can averagely achieve the highest signal-to-noise ratio (PSNR) and the structural similarity index (SSIM).

In Figure 2, an example of the denoised results for the compared algorithms are given. From Figure 2, we can see that our algorithm and the algorithms [10, 13] have similar visual effects, while the image restored by the algorithm [9] is somewhat blur.

To evaluate the performance of the compared algorithms quantitatively, the indexes PSNR and SSIM are used to measure the similarity between the denoised image and the original noisy-free image. The resulting values are displayed in Table 1. Since the code provided by the authors of [10] can only process the images of same width and height, some of the PSNR and SSIM of the algorithm [10] can not be listed in Table 1. From Table 1, we can see that our algorithm and the algorithm [13] have almost same PSNR and SSIM values averagely, which outperform ones of the algorithms in [9] and [10]. But the runtime for our algorithm is almost one-fifth of the runtime for the algorithm [13].
\begin{figure}[!t]
	\centering
	\subfloat[]{\includegraphics[scale=0.08]{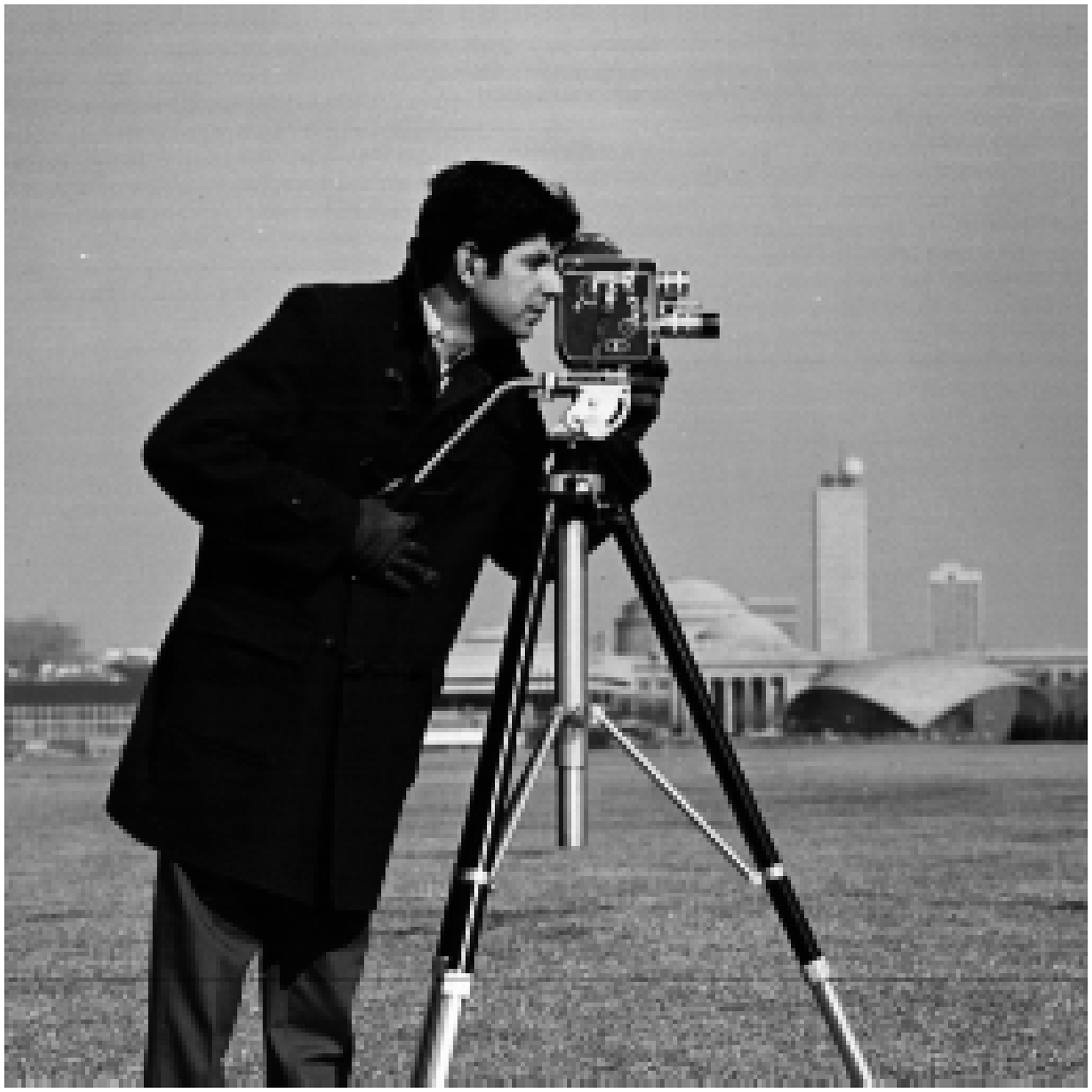}%
		\label{1}}
	\hfil
	\subfloat[]{\includegraphics[scale=0.08]{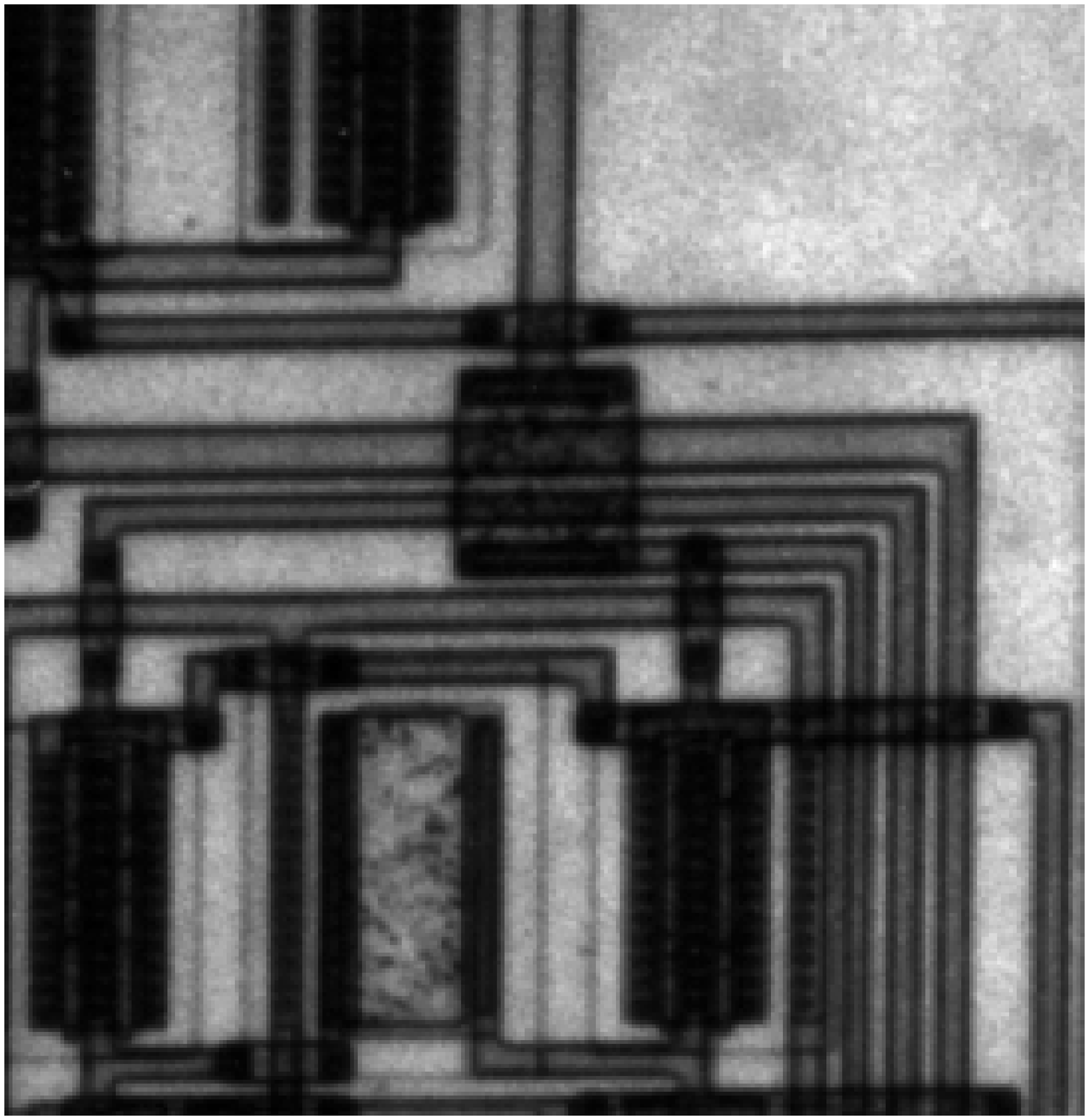}%
		\label{2}}
	\hfil
	\subfloat[]{\includegraphics[scale=0.08]{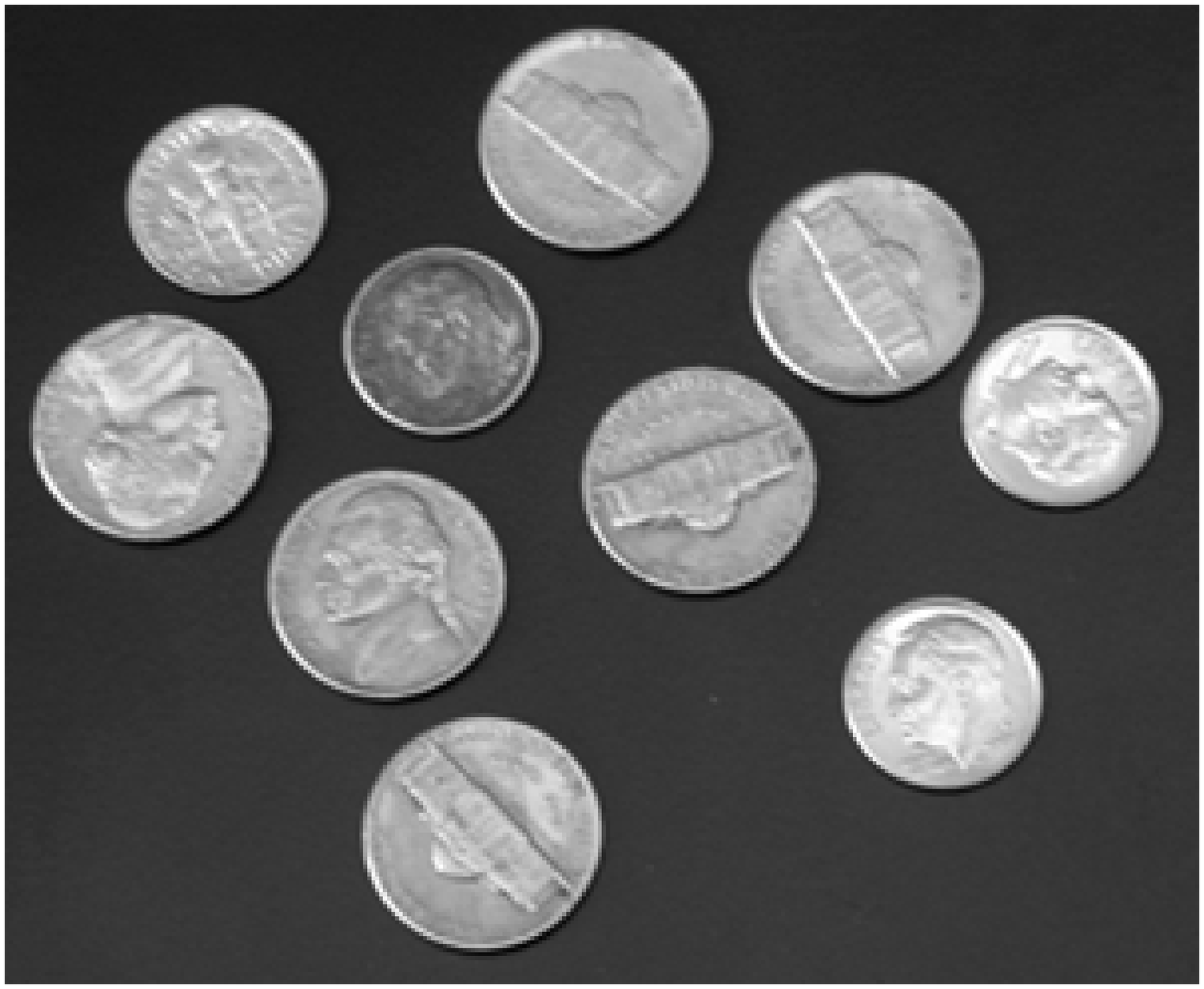}%
		\label{3}}
	\hfil
	\subfloat[]{\includegraphics[scale=0.08]{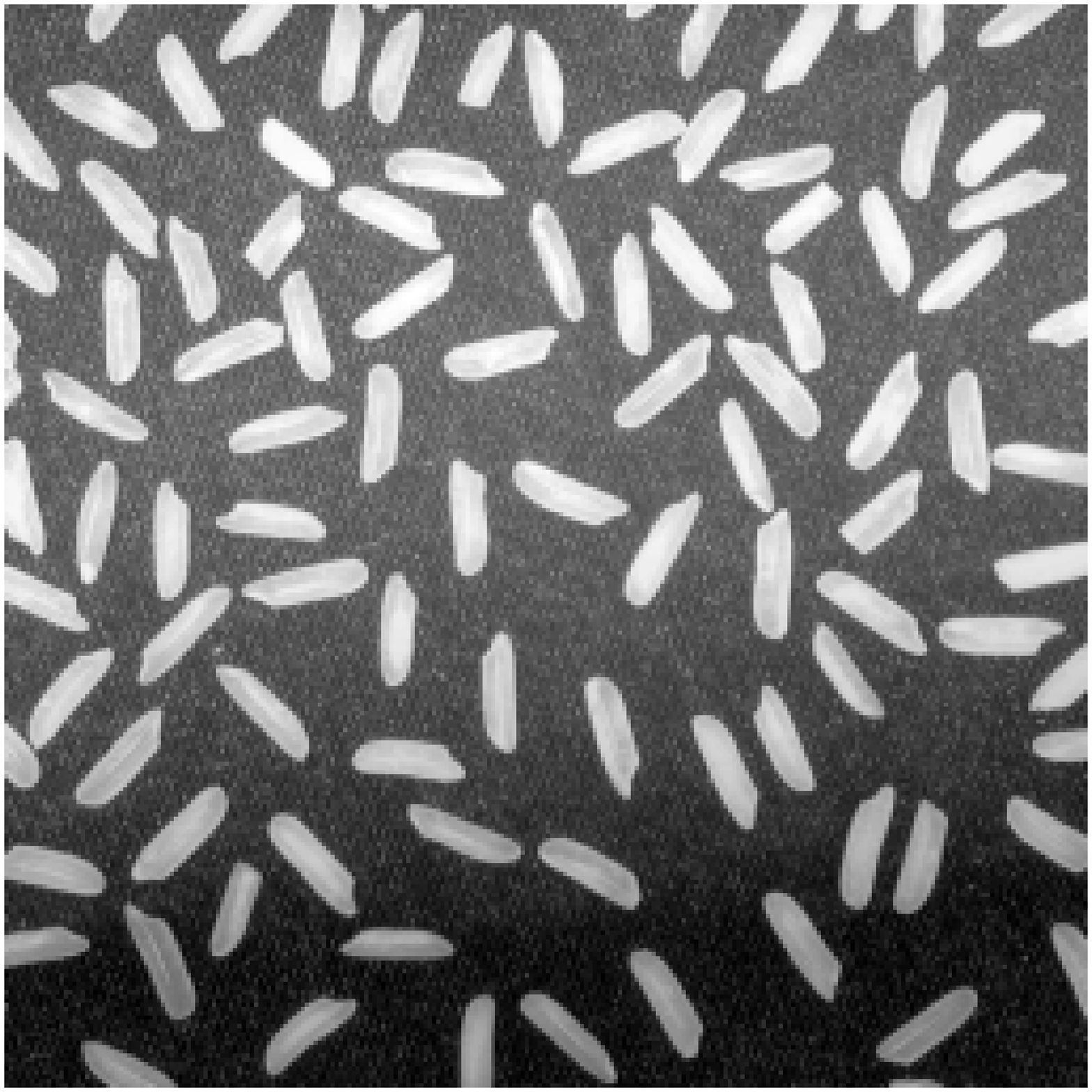}%
		\label{4}}
	\hfil
	\subfloat[]{\includegraphics[scale=0.08]{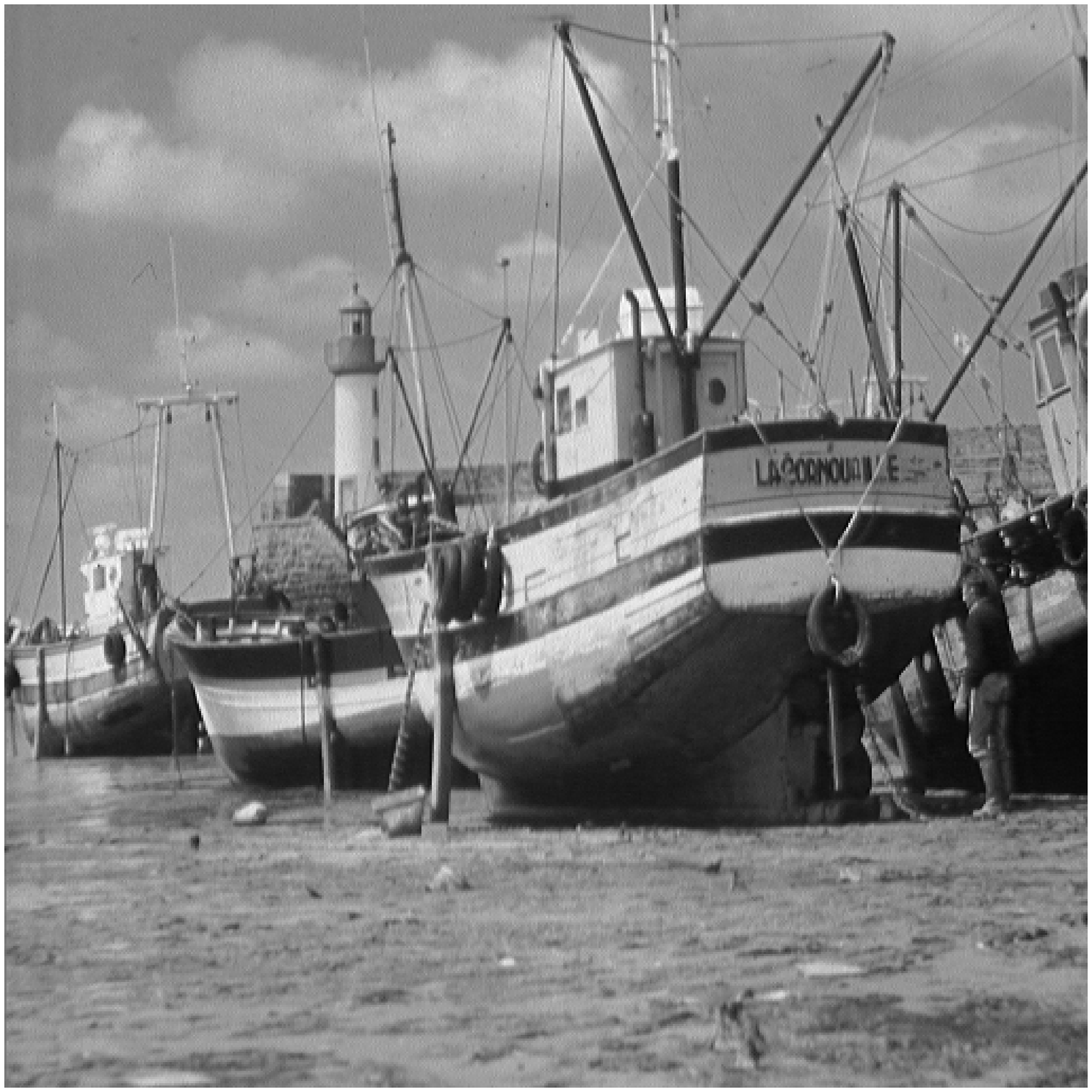}%
		\label{5}}\\
	
	\subfloat[]{\includegraphics[scale=0.08]{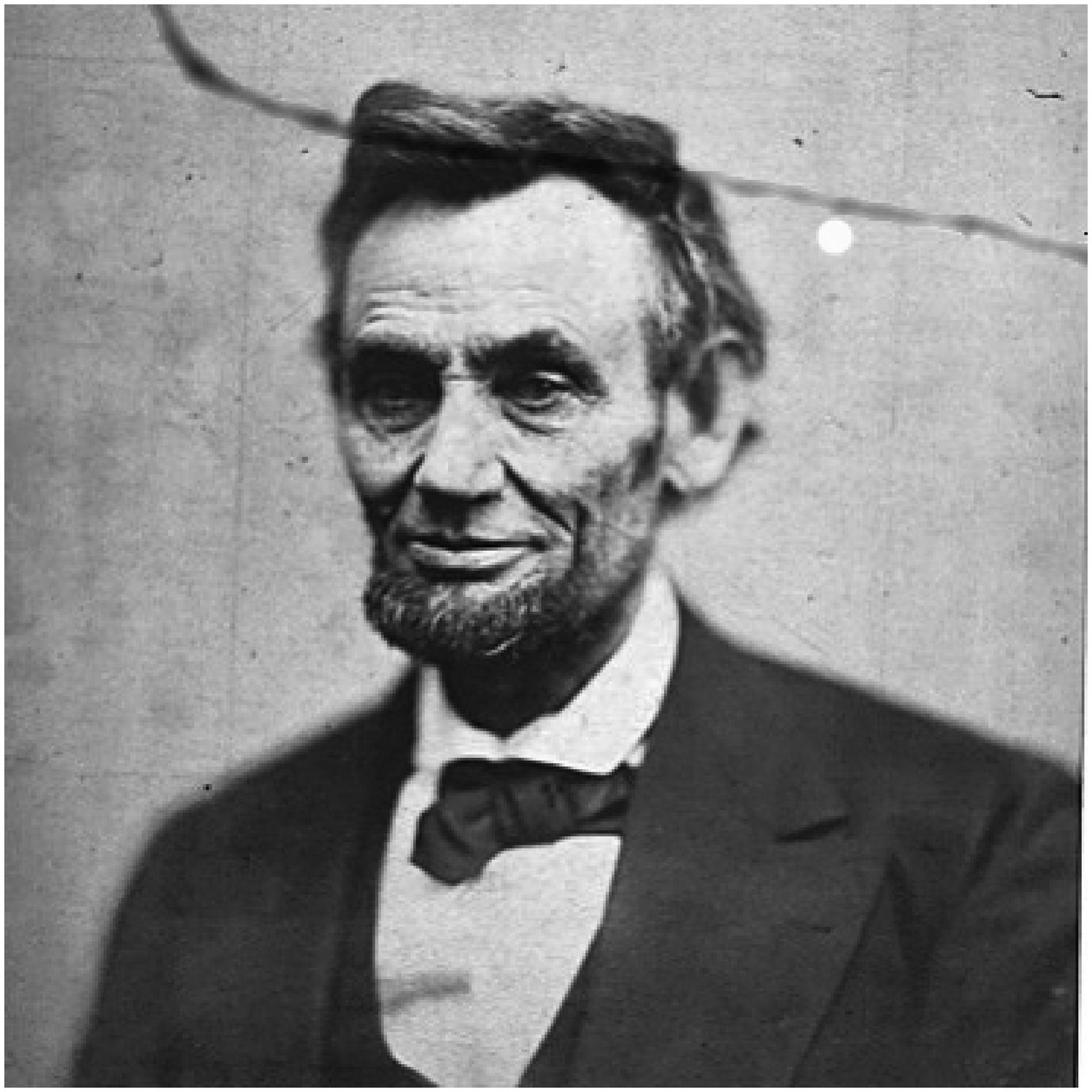}%
		\label{6}}
	\hfil
	\subfloat[]{\includegraphics[scale=0.08]{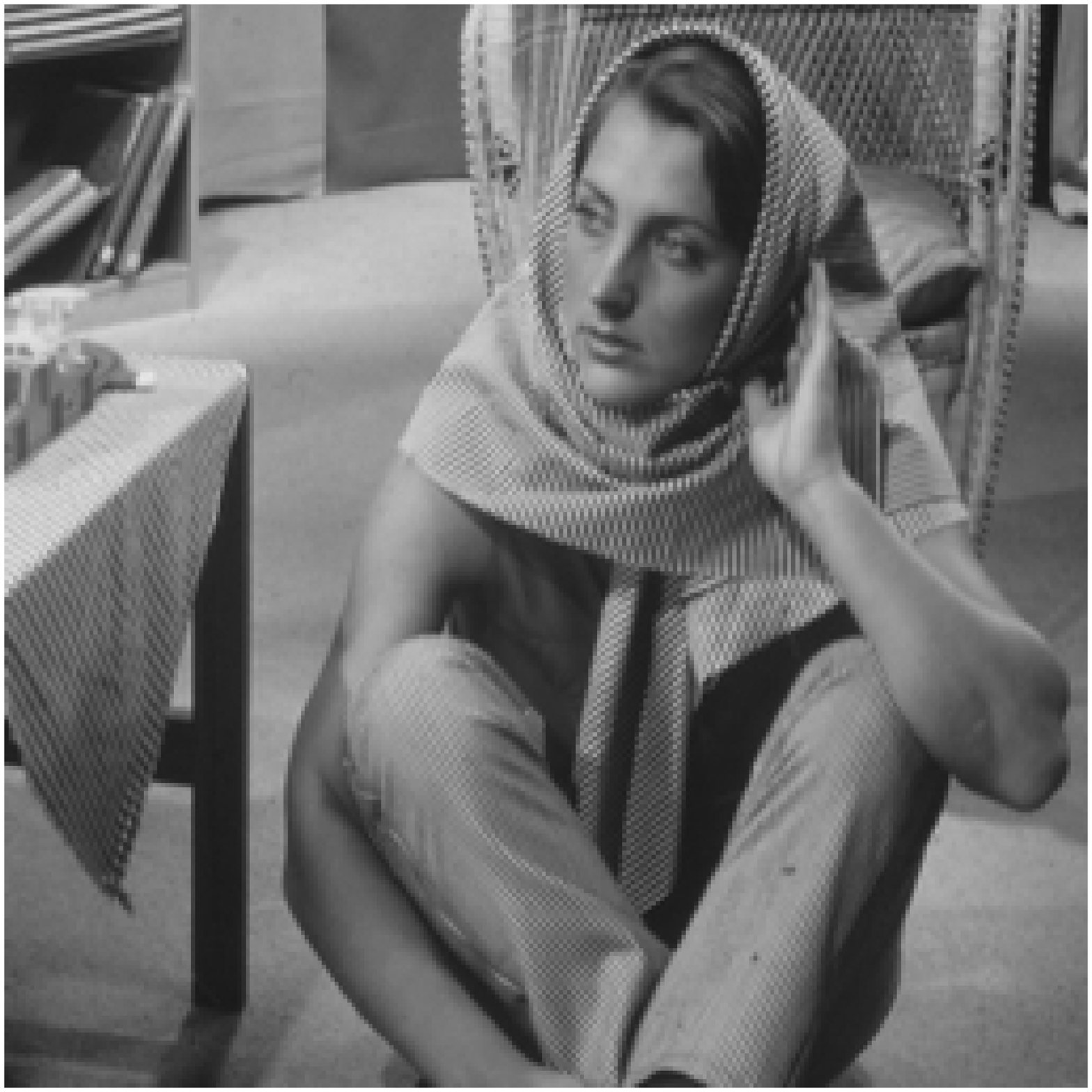}%
		\label{6}}
	\hfil
	\subfloat[]{\includegraphics[scale=0.08]{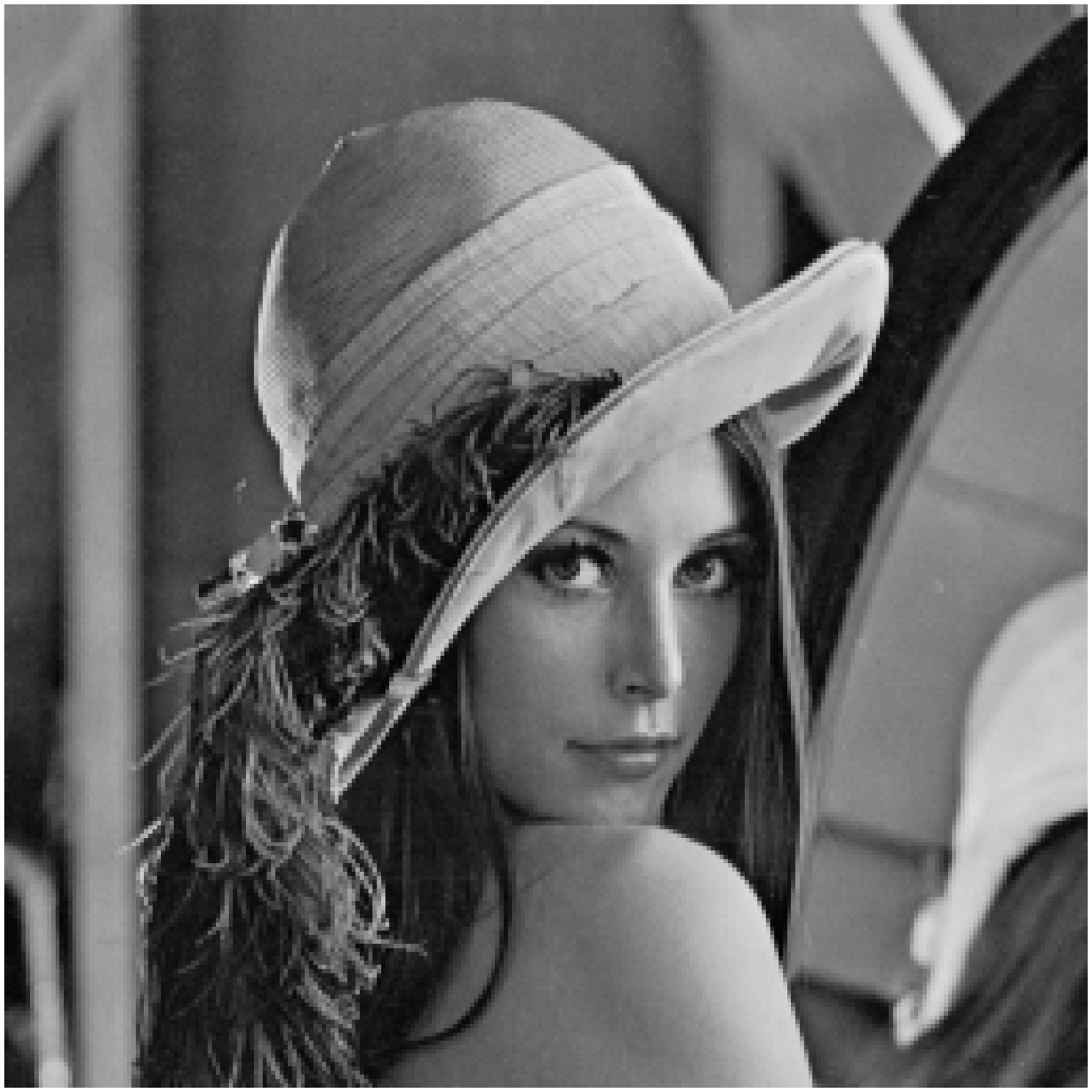}%
		\label{6}}
	\hfil
	\subfloat[]{\includegraphics[scale=0.08]{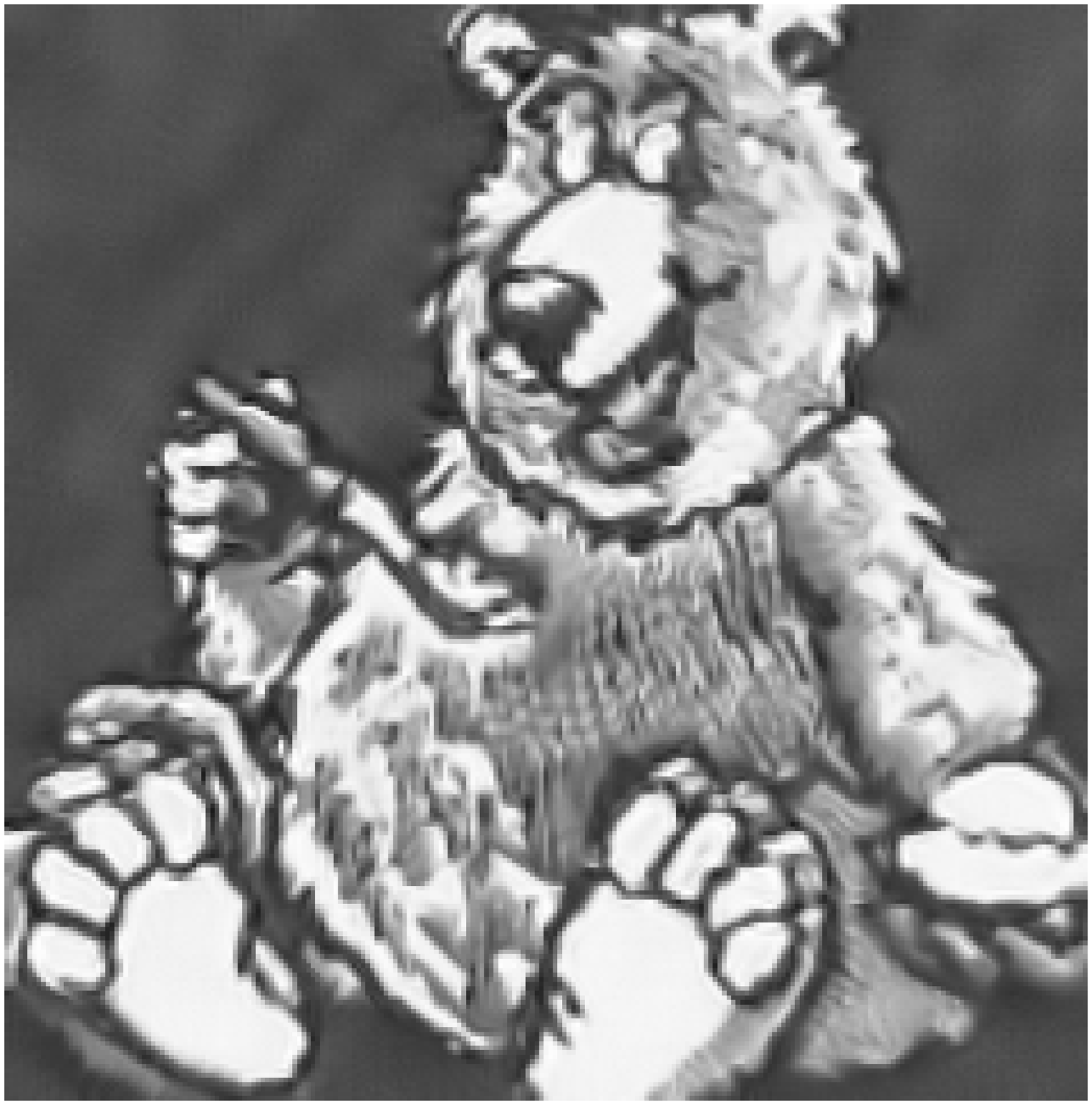}%
		\label{6}}
	\hfil
	\subfloat[]{\includegraphics[scale=0.08]{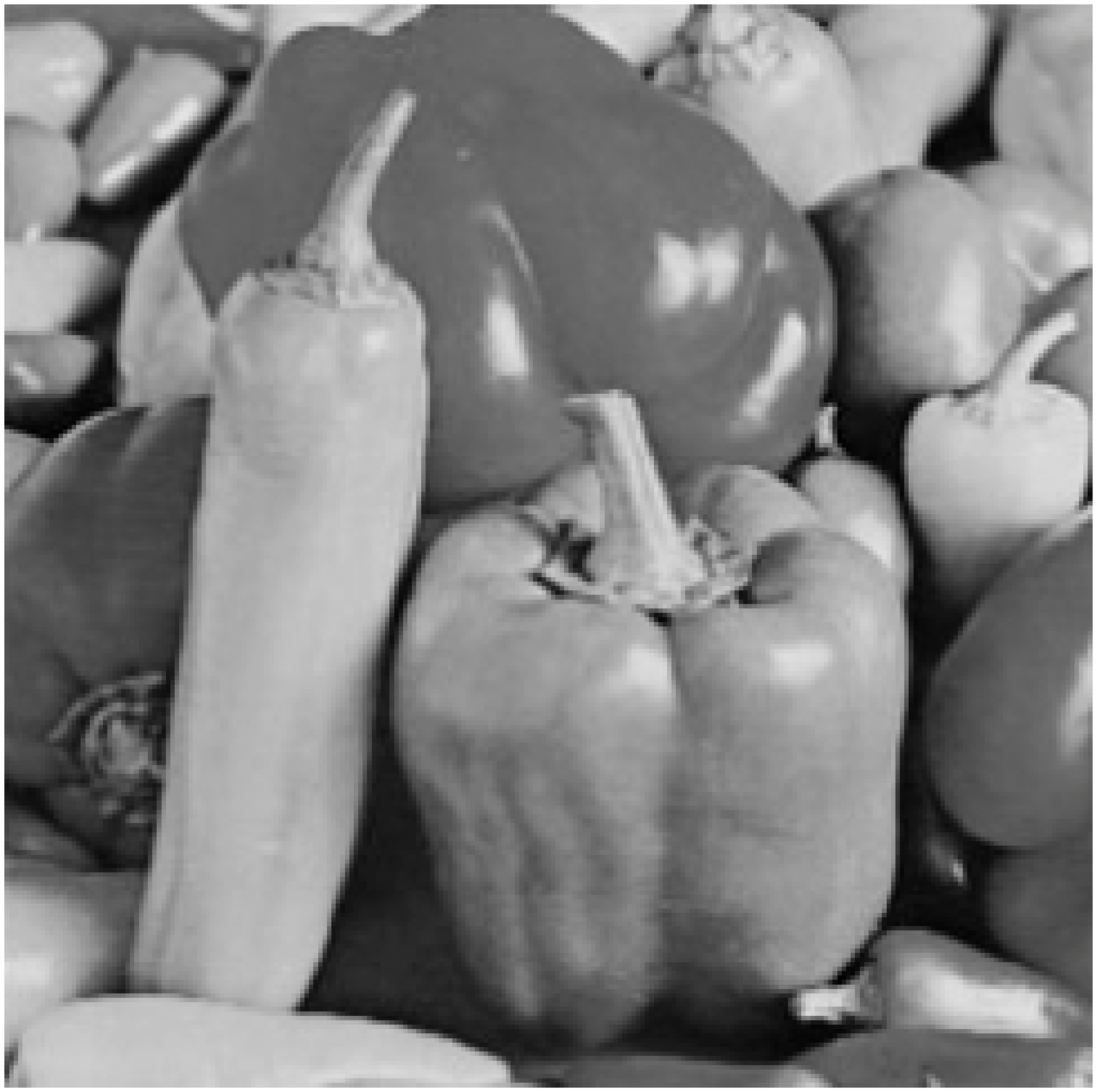}%
		\label{6}}
	\caption{original images.}	
\end{figure}

\begin{figure}[!t]
	\centering
	\subfloat[Noisy image]{\includegraphics[scale=0.2]{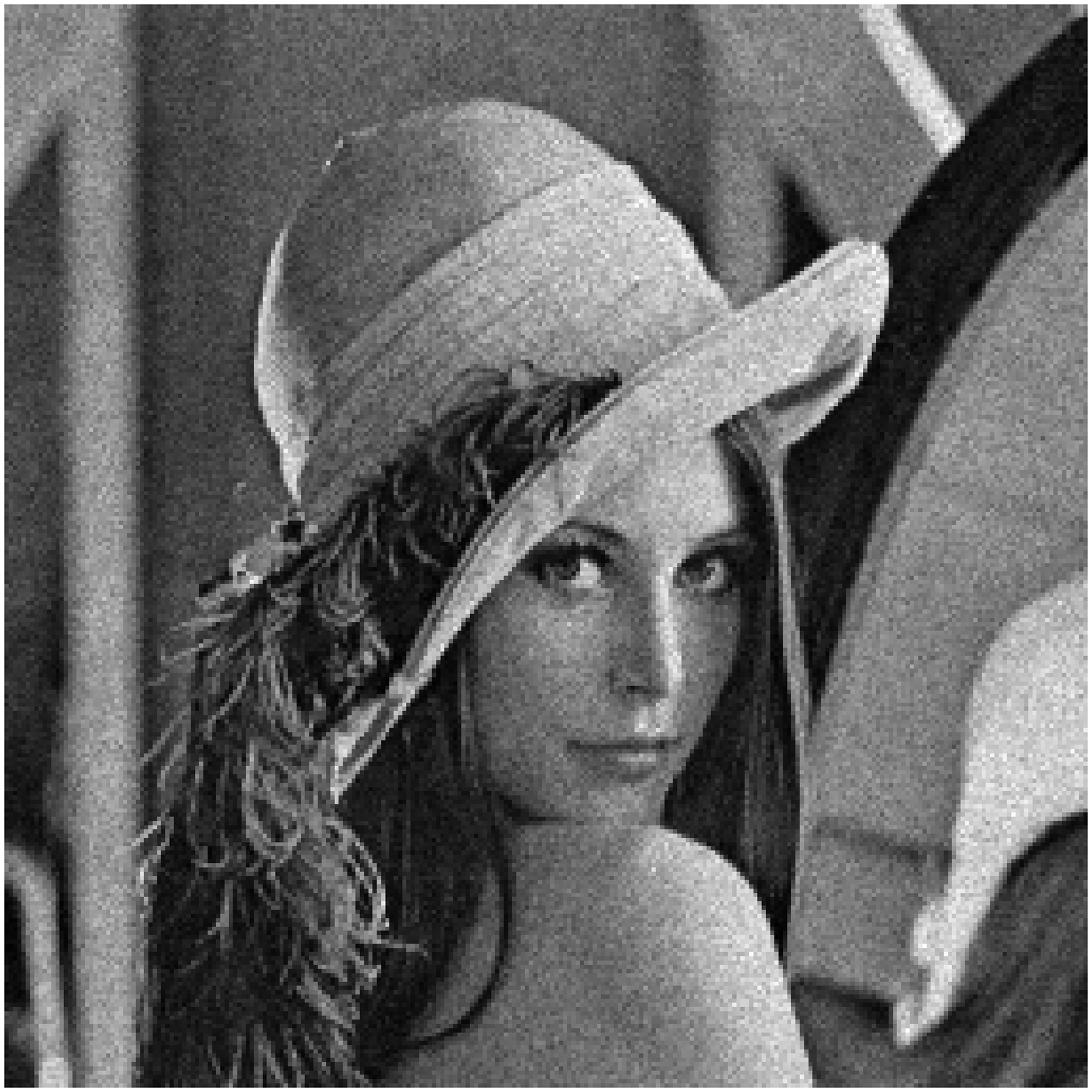}%
		\label{1}}
	\\
	\subfloat[Scheme in \text{[9]}]{\includegraphics[scale=0.2]{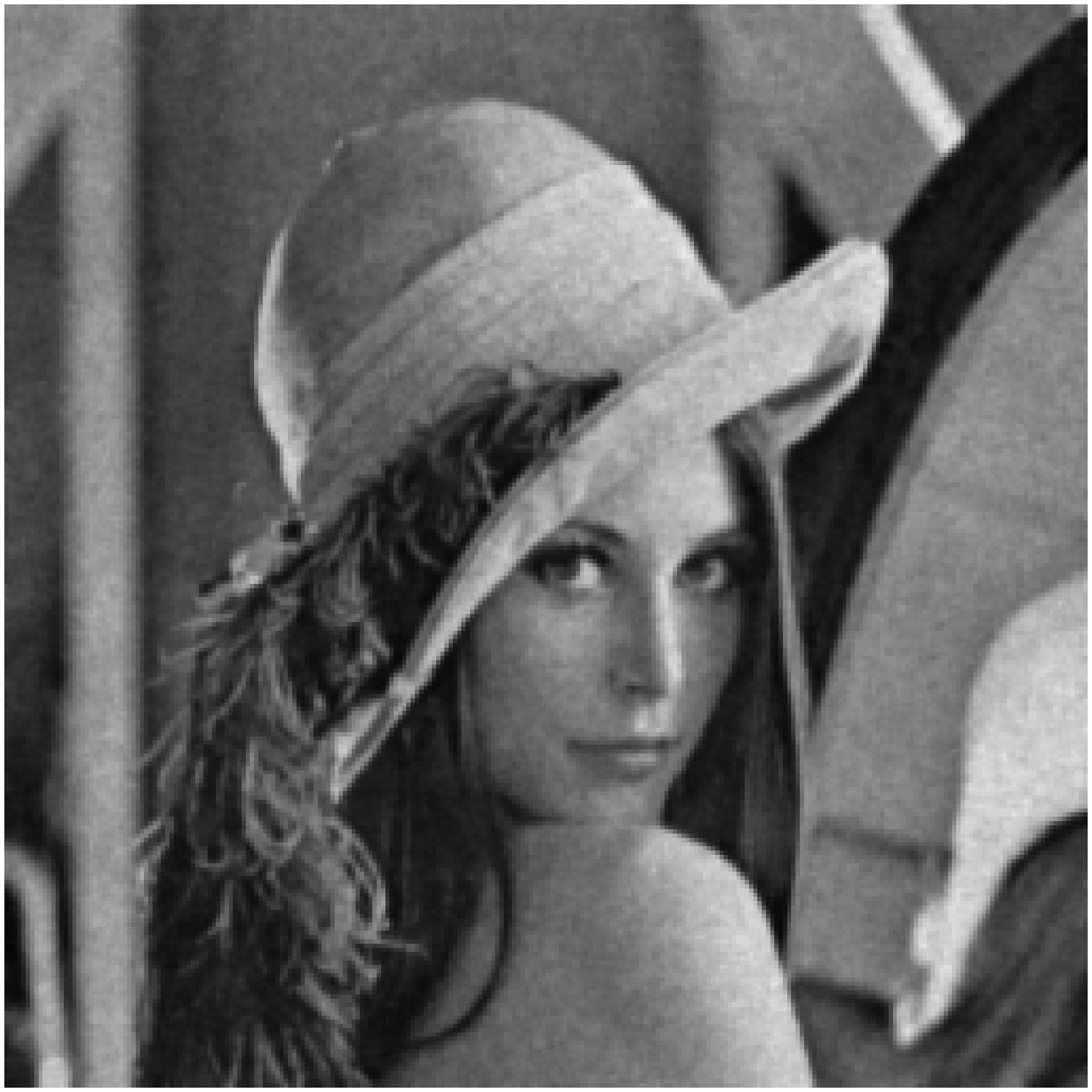}%
		\label{2}}
	\hfil
	\subfloat[Scheme in \text{[10]}]{\includegraphics[scale=0.2]{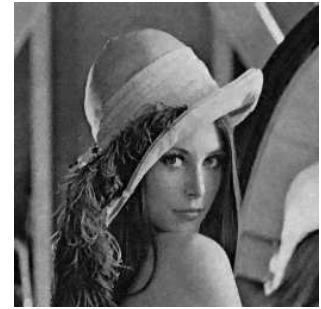}%
		\label{3}}
	\hfil\\
	\subfloat[Scheme in \text{[13]} ]{\includegraphics[scale=0.2]{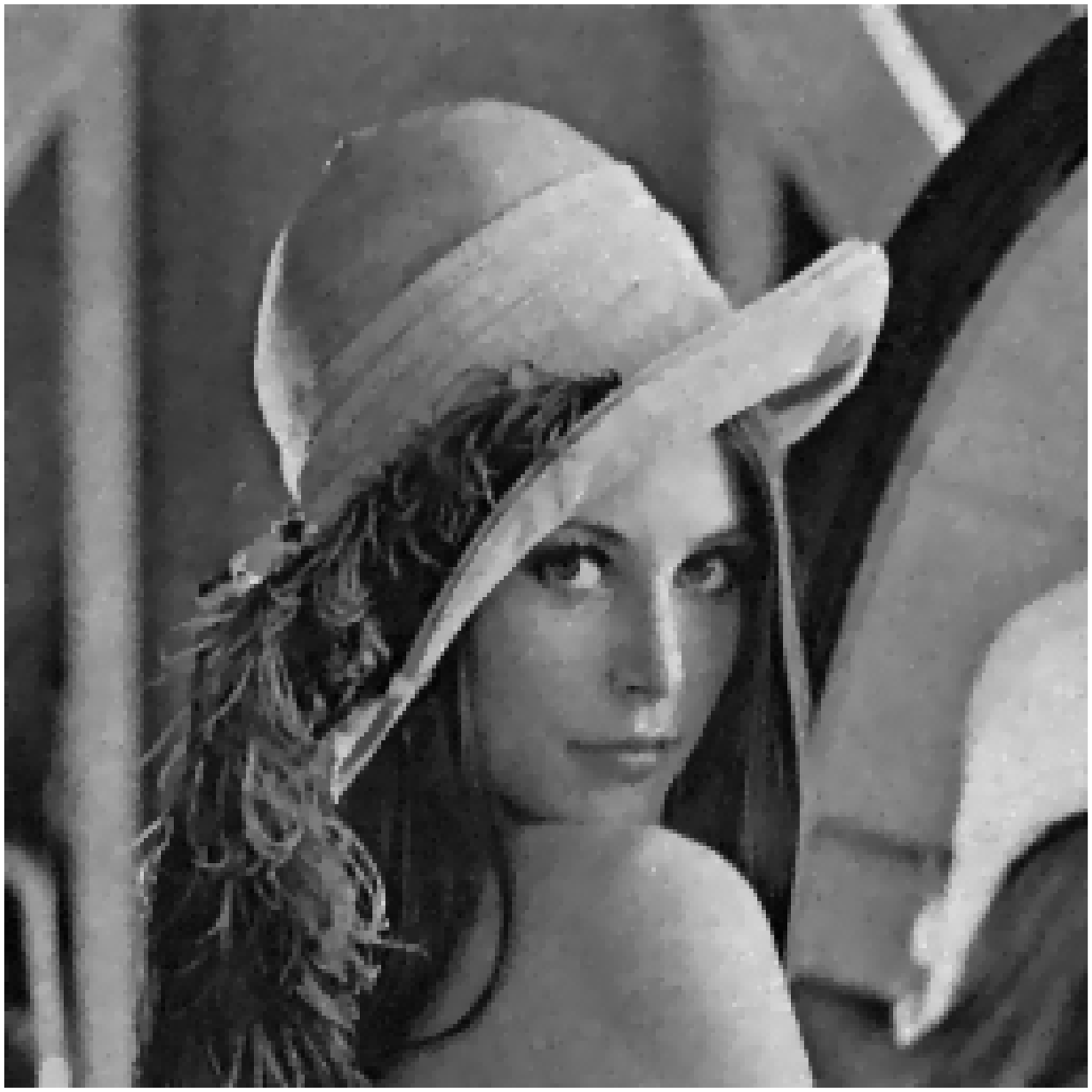}%
		\label{4}}
	\hfil
	\subfloat[Ours ]{\includegraphics[scale=0.2]{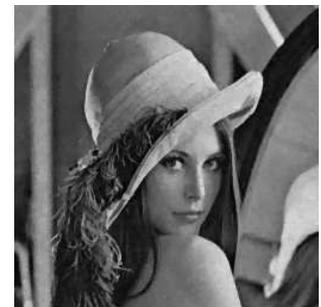}%
		\label{5}}
	
	\caption{Simulation results for the four schemes.}
	
\end{figure}

\begin{table*}[htbp]
	\centering
	\caption{The PSNR and SSIM of the compared four schemes}
	\centering
	\begin{tabular}{cc|cc|cc|cc|cc}
		\hline
		\multicolumn{2}{c}{Image} & \multicolumn{2}{c}{[9]} & \multicolumn{2}{c}{[10]} & \multicolumn{2}{c}{[13]} & \multicolumn{2}{c}{Ours} \\
		\hline
		ID    & Size  & PSNR  & SSIM  & PSNR  & SSIM  & PSNR  & SSIM  & PSNR  & SSIM \\ \hline
		a     & 256*256 & 28.04 & 0.8   & 32.21 & 0.87  & 32.18 & 0.88  & 32.17 & 0.88 \\
		b     & 272*280 & 33    & 0.92  &       &       & 32.85 & 0.91  & 32.86 & 0.91 \\
		c     & 300*246 & 31.27 & 0.88  &       &       & 33.50 & 0.92  & 33.51 & 0.92 \\
		d     & 256*256 & 27.71 & 0.85  & 30.61 & 0.9   & 30.34 & 0.87  & 30.36 & 0.87 \\
		e     & 512*512 & 30.58 & 0.82  & 31.76 & 0.83  & 31.94 & 0.84  & 31.94 & 0.84 \\
		f     & 371*371 & 31.89 & 0.82  & 32.83 & 0.82  & 33.23 & 0.84  & 33.21 & 0.84 \\
		g     & 256*256 & 30.41 & 0.84  & 23.31 & 0.82  & 31.59 & 0.88  & 31.60 & 0.88 \\
		h     & 256*256 & 30.55 & 0.86  & 32.41 & 0.88  & 33.06 & 0.90  & 33.08 & 0.90 \\
		i     & 237*240 & 28.64 & 0.89  &       &       & 30.89 & 0.94  & 31.00 & 0.93 \\
		j    & 243*242 & 32.45 & 0.86  &       &       & 33.84 & 0.91  & 33.75 & 0.91 \\ \hline
		\multicolumn{2}{r}{Average} \vline & 30.46 & 0.85  &  30.52     &  0.85     & 32.34 & 0.89  & 32.35 & 0.89 \\ \hline
		\multicolumn{2}{r}{Average time (s)} & \multicolumn{2}{r}{2.5703} & \multicolumn{2}{r}{66.2698} & \multicolumn{2}{r}{2.1328} & \multicolumn{2}{r}{0.4515} \\
		\hline
	\end{tabular}%
\end{table*}%

\section{Conclusion}
In this paper, we proved the existence and uniqueness of the LCA model by a new method. The semi-implicit scheme was designed to discretize the gradient flow, which can guarabtee the restored image to be positive in the image domain and allows for a lage time step. Experiments show that our method can numerically solve the LCA model  quickly and effectively.


%

\section*{Acknowledgment}

The authors would like to thank...

\ifCLASSOPTIONcaptionsoff
  \newpage
\fi


\begin{thebibliography}{1}

\bibitem{1}
S. W. Hasinoff, \emph{Photon, Poisson Noise}, Computer Vision, pp 608-610.
\bibitem{2}
B. Zhang, J. M. Fadili, and J.-L. Starck, \emph{Wavelets, ridgelets, and curvelets for Poisson noise removal}, IEEE TRANSACTIONS ON IMAGE PROCESSING, vol. 17, no. 7, pp. 1093-1108, Jul. 2008.
\bibitem{3}
C.-A. Deledalle, F. Tupin, and L. Denis, \emph{Poisson NL means: unsupervised non local means for Poisson noise}, in 2010 IEEE 17th Int. Conf. on Image Processing, Hong Kong, pp. 801-804, Sep. 2010.
\bibitem{4}
J. Salmon, \emph{On two parameters for denoising with non-local means}, IEEE SIGNAL PROCESSING LETTERS, vol. 17, no. 3, Mar. 2010.
\bibitem{5}
J. Salmon, C.-A. Deledalle, R.Willett, and Z. Harmany, \emph{Poisson noise reduction with non-local PCA}, in IEEE Int. Conf. Acoust., Speech and Signal Processing (ICASSP), 2012, pp. 1109–1112.
\bibitem{6}
André A. B. and Nelson D. A. M., \emph{A Nonlocal Poisson Denoising Algorithm Based on Stochastic Distances}, IEEE SIGNAL PROCESSING LETTERS, vol. 20, no. 11, pp. 1010-1013, Nov. 2013.
\bibitem{7}
M. Mäkitalo and A. Foi, \emph{Optimal inversion of the Anscombe transformation in low-count Poisson image denoising}, IEEE TRANSACTIONS ON IMAGE PROCESSING, vol. 20, no. 1, pp. 99–109, Jan. 2011.
\bibitem{8}
M. Mäkitalo and A. Foi, \emph{Optimal Inversion of the Generalized Anscombe Transformation for Poisson-Gaussian Noise}, IEEE TRANSACTIONS ON IMAGE PROCESSING, vol. 22, no. 1, pp. 91-103, Jan. 2013.
\bibitem{9} Le, T.; Chartrand, R.; Asaki, T. J., \emph{A Variational Approach to Reconstructing Images Corrupted by Poisson Noise}, JOURNAL OF MATHEMATICAL IMAGING AND VISION, Vol. 27, Iss. 3, pp. 257-263, Apr. 2007.
\bibitem{10} Chan, R. H. and Chen, Ke, \emph{Multilevel algorithm for a Poisson noise removal model with total-variation regularization}, INTERNATIONAL JOURNAL OF COMPUTER MATHEMATICS, vol. 84, iss. 8, pp. 1183-1198, 2007.
\bibitem{11} Zanella, R.; Boccacci, P.; Zanni, L.; et al., \emph{Efficient gradient projection methods for edge-preserving removal of Poisson noise}, INVERSE PROBLEMS, vol. 25, iss. 4, Arp. 2009.
\bibitem{12} Setzer, S., Steidl, G., Teuber, T., \emph{Deblurring Poissonian images by split Bregman techniques}, Journal of Visual Communication and Image Representation, vol. 21, pp. 193–199, 2010.
\bibitem{13} Figueiredo, M. A. T. and Bioucas-Dias, J. M., \emph{Restoration of Poissonian Images Using Alternating Direction Optimization}, IEEE TRANSACTIONS ON IMAGE PROCESSING, Vol. 19, Iss. 12, Pp. 3133-3145, Dec. 2010.
\bibitem{14} Liu, X. and Huang, L., \emph{Total bounded variation-based Poissonian images recovery by split Bregman iteration}, MATHEMATICAL METHODS IN THE APPLIED SCIENCES, vol. 35, iss. 5, pp. 520-529, Mar. 2012.
\bibitem{15} Wang, X.; Feng, X.; Wang, W.; et al., \emph{Iterative reweighted total generalized variation based Poisson noise removal model}, APPLIED MATHEMATICS AND COMPUTATION, vol. 223, pp. 264-277, Oct. 2013.
\bibitem{16} W. Zhou and Q. Li, \emph{Adaptive total variation regularization based scheme for Poisson noise removal}, MATHEMATICAL METHODS IN THE APPLIED SCIENCES, vol. 36, iss. 3, pp. 290-299, Feb. 2013.
\bibitem{17} R. Chan, H. Yang, and T. Zeng, \emph{A Two-Stage Image Segmentation Method for Blurry Images with Poisson or Multiplicative Gamma Noise}, SIAM Journal on Imaging Sciences, vol. 7, no. 1, pp. 98–127, 2014.
\bibitem{18} Jiang, L.; Huang, J.; Lv, X. et al., \emph{Alternating direction method for the high-order total variation-based Poisson noise removal problem}, NUMERICAL ALGORITHMS, vol. 69, iss. 3, pp 495-516, Jul. 2015.
\bibitem{19} Karimi, D. and Ward, R., \emph{A denoising algorithm for projection measurements in cone-beam computed tomography}, COMPUTERS IN BIOLOGY AND MEDICINE, vol. 69, pp. 71-82, Feb. 2016.
\bibitem{20} Mansouri, A.; Deger, F.; Pedersen, M.; et al, \emph{An adaptive spatial-spectral total variation approach for Poisson noise removal in hyperspectral images}, SIGNAL IMAGE AND VIDEO PROCESSING, vol. 10, iss. 3, pp. 447-454, Mar. 2016.
\bibitem{21} Liu, X., \emph{Augmented Lagrangian method for total generalized variation based Poissonian image restoration}, COMPUTERS and MATHEMATICS WITH APPLICATIONS, vol. 71, iss. 8, pp. 1694-1705, Apr. 2016.

\bibitem{22} L. Rudin, S. Osher, and E. Fatemi, \emph{Nonlinear total variation based noise removal algorithms,} Physica D, vol. 60, pp. 259-268, 1992.
\bibitem{23} Aubert, G.; Aujol, J., \emph{a variational approach to remove multiplicative noise}, SIAM JOURNAL ON APPLIED MATHEMATICS, vol. 68, iss. 4, pp. 925-946, 2008.

\end{thebibliography}
\end{document}